\documentclass[10pt]{article}

\usepackage{amsmath}
\usepackage{amssymb}
\usepackage{indentfirst}
\usepackage{graphics} 
\usepackage{color}
\usepackage[utf8]{inputenc}

\makeatletter

\@addtoreset{equation}{section}
\makeatother
\def\R{\mbox{\boldmath $R$}}

\def\<{\langle}
\def\>{\rangle}

\newtheorem{lem}{Lemma}[section]
\newtheorem{theo}{Theorem}[section]
\newtheorem{rem}{Remark}[section]
\newtheorem{pro}{Proposition}[section]

\makeatletter
   
   \@addtoreset{equation}{section}
\makeatother

\begin{document}
\title{\bf  Fast energy decay for damped wave equations\\ with a potential and rotational inertia terms}

\author{Ruy Coimbra Char\~ao\thanks{ruy.charao@ufsc.br}  \\{\small Department of Mathematics}, {\small Federal University of Santa Catarina} \\ {\small 88040-900, Florianopolis, SC,  Brazil,} 
\\
and\\Ryo Ikehata\thanks{ikehatar@hiroshima-u.ac.jp} \\ {\small Department of Mathematics}, \small{ Division of Educational Sciences}\\ {\small Graduate School of Humanities and Social Sciences} \\ {\small Hiroshima University} \\ {\small Higashi-Hiroshima 739-8524, Japan}}

\date{}
\maketitle
\begin{abstract}
We consider damped wave equations with a potential and rotational inertia terms. We study the Cauchy problem for this model in the one dimensional Euclidean space ${\bf R}$ and we obtain fast energy decay and $L^{2}$-decay of the solution itself as $t \to \infty$. Since we are considering this problem in the one dimensional space, we have no useful tools such as the Hardy and/or Poincar\'e inequalities. This causes significant difficulties to derive the decay property of the solution and the energy. A potential term will play a role for compensating these weak points. 
\end{abstract}
\section{Introduction}
\footnote[0]{Keywords and Phrases: Damped wave equation; potential; rotational inertia; one dimension; fast energy decay.}
\footnote[0]{2010 Mathematics Subject Classification. Primary 35L05; Secondary 35B40, 35C20, 35S05.}
We consider the following dissipative wave type equation that may model linear hydrodynamic problems
\begin{align}
& u_{tt} - u_{ttxx} - u_{xx} + V(x) u + u_{t}  = 0,\ \ \ (t,x)\in (0,\infty)\times {\bf R},\label{eqn}\\
& u(0,x)= u_0(x), \quad  u_{t}(0,x)= u_{1}(x),\ \ \ x\in{\bf R} ,\label{initial}
\end{align}
where $V(x) > 0$ is a potential specified later.
The initial data and the solution shall be taken as a real value.

\noindent
{\bf Notation.} $f \in {\rm BC}({\bf R})$ implies $f(x)$ is continuous and bounded in ${\bf R}$. $f \in {\rm BC}^{1}({\bf R})$ implies $f \in {\rm BC}({\bf R})$, and $f' \in {\rm BC}({\bf R})$, and $f \in {\rm BC}^{2}({\bf R})$ implies $f \in {\rm BC}^{1}({\bf R})$, and $f'' \in {\rm BC}({\bf R})$.To simplify the notation in the calculations, the following symbols are sometimes used with dimensional generalization in mind: $\nabla = \frac{\partial}{\partial x}$, $\Delta = \frac{\partial^{2}}{\partial x^{2}}$, $\Delta^2 = \frac{\partial^{4}}{\partial x^{4}}$, and $\Delta u_x = \frac{\partial^{3}u}{\partial x^{3}}$. Moreover, we use the symbols $||\cdot||$ and $(\cdot,\cdot) = (\cdot,\cdot)_{L^{2}}$ to denote the norm and  the inner product in $L^2=L^2({\bf R})$, respectively. By $\Vert \cdot\Vert_{H^{l}}$ we denote the usual $H^{l}({\bf R})$-norm. Finally, $\Vert f\Vert_{\infty}$ means the usual $L^{\infty}$-norm of $f \in L^{\infty}({\bf R})$.
\\

Concerning the existence of a unique solution to problem \eqref{eqn}-\eqref{initial}, by a similar argument to \cite[Proposition 2.1]{ITY} based on the Lumer-Phillips Theorem one can find that the problem (1.1)-(1.2) with initial data $(u_0, u_1) \in H^{2}({\bf R}) \times H^{2}({\bf R}) $  has a unique strong solution
\[u \in C^{2}([0,\infty);H^{2}({\bf R}))\]
satisfying
\begin{equation}\label{1-1}
E_{u}(t) + \int_{0}^{t}\Vert u_{s}(s,\cdot)\Vert^{2}ds = E_{u}(0),
\end{equation}
where
\[E_{u}(t)=E(t)=:\frac{1}{2}\left(\Vert u_{t}(t,\cdot)\Vert^{2} + \Vert \nabla u_{t}(t,\cdot)\Vert^{2} + \Vert\nabla u(t,\cdot)\Vert^{2} + \Vert\sqrt{V(\cdot)}u(t,\cdot)\Vert^{2}\right)\] 
represents the total energy with respect to the equation \eqref{eqn} (see Appendix below (cf. \cite[Th\'eor$\grave e$me VII.5]{B} and \cite{Ikawa, LI, pazy, WC-2})). 
Furthermore, one can get the regularity property such that for $[u_{0},u_{1}] \in H^{4}({\bf R})\times H^{4}({\bf R})$, the corresponding solution $u(t,x)$ satisfies $u \in C^{2}([0,\infty);H^{4}({\bf R}))$, which guarantees the validity of all computations developed in this paper via density argument (for example, Lemma \ref{lem34} below).
\\

First, let us discuss the relevant prior results and the motivation for the study.\\ 
Equation \eqref{eqn} can also be considered as an equation based on the generalized IBq equation with a certain friction term $u_{t}$, as studied, for example, in reference \cite{WX} and \cite{WLZ}, plus a potential term $V(x)u$ that depends on spatial variables. Of course, it should be pointed out that the original IBq equation is considered with nonlinear terms as follows. For $\nu > 0$, 
\[u_{tt} - \Delta u_{tt} - \Delta u + \nu u_{t}  = \Delta f(u), f(u_{t}), \Delta f(u, u_t).\]
The linearized version of this equation without the nonlinear term $\Delta f(u)$ and/or $f(u_{t})$ will be deceased here, but of course, these studies should be useful enough to study the equation with the added nonlinear terms. Incidentally, for more information on the historical topics and physical meaning of these Bq or IBq equations, refer to \cite{WC, WX, XW} and the references therein. In several articles authors studied asymptotic $L^2$-norm decay of solutions with hydrodynamical dissipation or Stokes damped term e.g. \cite{WX-2, WX}, but in the article \cite{WANG} under some assumptions it is proved that the $L^{\infty}$--norm of small solution of the Cauchy problem for the Boussinesq equation  decays to zero as $t$ tends to infinity. In \cite{PE} the existence and blow-up of solution of Cauchy problem are considered for the generalized damped multidimensional Boussinesq equation.
\\
In particular, about the scattering and the growth properties of the solution of the (non-damped) IMBq equation:
\[u_{tt}-\Delta u_{tt} -\Delta u = \Delta f(u),\]
one can refer to the papers \cite{WC-2} for $f(u) \ne 0$ and \cite{LI} for $f(u) = 0$. For a related optimal decay estimates for the solutions of the dissipative "plate" type equation including the linear case:
\[u_{tt}-\Delta u_{tt} + \Delta^2 u + u_{t} = 0,\]
one can refer the reader to \cite{LK, LC}, and the references therein.

Now, according to the results of \cite{WX} and \cite{WLZ}, it is pointed out that the so-called regularity-loss structure is inherent in the decay structure of the solution of the IBq equation without this nonlinear term. In fact, in \cite{WX-2} the authors observe that the dissipative structure via hydrodynamical
term  of the linearized IBq equation is of the regularity-loss type.
Moreover, in the article \cite{WZ} the authors show that the solution of generalized  Boussinesq equation for dimension $n \geq 2$ can be approximated by the linear solution as time tends to infinity. So, it would be interesting to see the effect of adding a potential term $V(x)u$ on such a structure. In \cite{WX}, the authors are dealing with constant coefficient case, so analysis through Fourier transform is possible, but here, due to the presence of spatial variable-dependent potential terms, we need to develop a different analysis method, which will make it more difficult. Incidentally, in this paper, we deal with the one dimension, but of course we will perform calculations that are valid in general dimensions as well. However, the real new decay estimate formulas, in which the influence of the potential term is pronounced, are limited to the $1D$ case. It should be pointed out that the regularity-loss structure itself has been first discovered and named by Professor Shuichi Kawashima in the papers \cite{HK, IHK}.

In addition,  we note again that the following linear equation under effects of a disturbance of the physical environment by a potential and a Stokes damped term can be considered to study hydrodynamical models in place of the the nonlinear model considered by Boussinesq: 
 $$  u_{tt} - u_{ttxx} - u_{xx} + V(x) u + u_{t}  =0,\ \ \ (t,x)\in (0,\infty)\times {\bf R}.$$ \\
It is important to note that the first attempt of mathematical modeling of shallow water was by  Joseph Valentin Boussinesq in 1871, who took into account the vertical accelerations of the fluid particles and admitted a solitary solution, but he did not consider terms of derivative products considering only nonlinear
first-order effects. These equations are called classical Boussinesq
equations. In order to obtain better models for the hydrodynamic phenomena, other equations were
developed, considering higher order terms: class Boussinesq equations.
For example equations of Korteweg-de Vries, Green-Naghdi [J. Mech
Fluid. 1976], Peregrine [J. Fluid Mech 1967 - long waves] and F. Serre
\cite{Serre} (see also \cite{CH}). Other equations related to Boussinesq are the so-called IBq and
IMBq (see references in these articles mentioned above).
\\
\vspace{0.2cm}
Before starting our main result, we shall mention the assumption on the potential $V(x)$, which can be admitted in our argument.\\
On the potential $V(x)$, we assume that\\
{\bf (V.1)}\,\, $V \in {\rm BC}^{2}({\bf R})$, $V(x) > 0$ and $\vert V'(x)\vert \leq \alpha V(x)$ ($\forall x \in {\bf R}$) with some constant  $\alpha > 0$.\\
\underline{{\bf Example.}} One can show a typical example satisfying {\bf (V.1)}: $V(x) := V_{0}(1+x^{2})^{-\frac{\alpha}{2}}$ ($\alpha > 0$ and $V_{0} > 0$). However, one can not choose an exponential type such as $V(x) = e^{-x^{2}}$. It imposes some limits on the degree of decay in the spatial distance. \\
Our purpose is to find a fast decay rate of the total energy. Several experiences in analysing plate-type equations with rotational inertia terms is documented in \cite{CI, LC}, but careful consideration is needed here to compensate for the critical weakness of not being able to use Hardy-type inequalities. Our result reads as follows.
\begin{theo}\label{theorem-1}
Assume {\bf (V.1)} and $\alpha^{2}\Vert V\Vert_{\infty} < 1$. Let $u_{0} \in H^{3}({\bf R})$ and $u_{1} \in H^{3}({\bf R})$ satisfy
$$\frac{u_{0} +u_{1}-\Delta u_{1}}{\sqrt{V(\cdot)}} \in L^2({\bf R}).$$
Then, the strong solution $u(t,x)$ to problem \eqref{eqn}-\eqref{initial} satisfies
\[E_{u}(t) \leq C(\Vert u_{0}\Vert_{H^{3}}^{2} + \Vert u_{1}\Vert_{H^{3}}^{2} + \left\|\dfrac{u_{0}+u_{1}-\Delta u_{1}}{\sqrt{V(\cdot)}}\right\|^{2})(1+t)^{-2},\quad (t \geq 0),\]
\[\Vert u(t,\cdot)\Vert^{2} \leq C(\Vert u_{0}\Vert_{H^{3}}^{2} + \Vert u_{1}\Vert_{H^{3}}^{2} + \left\|\dfrac{u_{0}+u_{1}-\Delta u_{1}}{\sqrt{V(\cdot)}}\right\|^{2})(1+t)^{-1},\quad (t \geq 0),\]
where $C = C(\alpha, \Vert V\Vert_{\infty}, \Vert V'\Vert_{\infty}, \Vert V''\Vert_{\infty}) > 0$ is a constant depending only on the quantities $\alpha$, $\Vert V\Vert_{\infty}$, $\Vert V'\Vert_{\infty}$ and $\Vert V''\Vert_{\infty}$.
\end{theo}
The proof of this theorem will be done in the next sections.
\begin{rem}{\rm It should be pointed out that even when the initial value is smooth enough, e.g., $(u_{0},u_{1}) \in C_{0}^{\infty}({\bf R})\times C_{0}^{\infty}({\bf R})$, the result of the theorem is new enough. In this case the assumption $\frac{u_{0} +u_{1}-\Delta u_{1}}{\sqrt{V(\cdot)}} \in L^2({\bf R})$ automatically holds.}
\end{rem}
\begin{rem}{\rm Comparing the results of Theorem \ref{theorem-1} above and \cite[Theorem 4.1]{WLZ}, both the energy and the $L^2$ norm of the solution derive faster decay rates in Theorem \ref{theorem-1}. This is still an effect of the addition of the potential term $V(x)u$. Here we are limited to one spatial dimension, but a slightly different method would be needed to generalize this to a more general dimension.}
\end{rem}
\noindent
\underline{{\bf Example.}}\, For the concrete example $V(x) := V_{0}(1+x^{2})^{-\frac{\alpha}{2}}$, since $\Vert V\Vert_{\infty} = V_{0}$, the condition $\alpha^{2}\Vert V\Vert_{\infty} < 1$ implies $\alpha^{2}V_{0} < 1$. That is, if we choose large $\alpha \gg 1$, then we must choose small $V_{0}$, and vice versa. Incidentally, the large $\alpha > 0$ corresponds to the so-called short range potential.
\begin{rem}{\rm To get faster decay like $O(t^{-2})$ of the total energy one must assume further regularity on the initial data and the potential. This may reflect an essential effect from the rotational inertia term $-\Delta u_{tt}$. This can be seen by comparing it to the case of the damped wave equation with potential (see \cite[Theorem 1.3]{I-2}). In \cite[Theorem 1.3]{I-2}, the same decay order as in Theorem \ref{theorem-1} above can be derived under weaker regularity assumption $(u_{0}, u_{1}) \in H^{1}\times L^{2}$ to the one dimensional equation:
\[u_{tt} - \Delta u + V(x) u + u_{t} = 0.\]
}
\end{rem}
\begin{rem}{\rm Although Theorem \ref{theorem-1} definitely yields a fast decay rate, at this stage it is a future issue to be addressed as to the best decay rate. An interesting recent paper \cite{S} may be of interest in this regard. }
\end{rem}

The paper is organized as follows. In section 2, we will study slower decay of the total energy under weaker conditions on the initial data and the potential $V(x)$. Section 3 will be devoted to the proof of Theorem \ref{theorem-1} by using one particular multiplier method under high regularity condition on the initial data and potential. In Appendix we shall discuss the unique existence of solutions to problem \eqref{eqn}-\eqref{initial} based on the Lumer-Phillips theorem. \\


\section{Main Estimates}
We first prepare important lemma to proceed our argument in this paper. A main idea comes from \cite[Lemma 3.4]{I-2}. 
\begin{lem}\label{lem21}\,Let $V \in {\rm BC}^{1}({\bf R})$ satisfies $V(x) > 0$ for $x\in{\bf R}$ and $[u_{0},u_{1}] \in H^{2}({\bf R})\times H^{2}({\bf R})$ satisfying
$$\frac{u_{0} +u_{1}- \Delta u_{1}}{\sqrt{V(\cdot)}} \in L^2({\bf R}).$$
Then the strong solution $u \in C^{2}([0,\infty);H^{2}({\bf R}))$ to problem \eqref{eqn}-\eqref{initial}
satisfies
\[\Vert u(t,\cdot)\Vert^{2} + \int_{0}^{t}\Vert u(s,\cdot)\Vert^{2}ds \leq J_{0}^{2}\quad (t \geq 0),\]
where $$J^{2}_{0}:= \dfrac{1}{2}\big(\|u_{0}\|^{2}+\|\nabla u_{0}\Vert^{2}\big)
 +\dfrac{1}{2}\left\|\dfrac{u_{0}+u_{1}-\Delta u_{1}}{\sqrt{V(\cdot)}}\right\|^{2}.$$
\end{lem} 
{\bf Proof:} 
Let us define $\displaystyle w(t,x) := \int_{0}^{t}u(s,x)\,ds$. Then $w_{t}(t,x) = u(t,x)$ and $w(0)=0$, $w_{t}(0)=u_{0}$ and the function $w(t,x)$ satisfies
\begin{align}
& w_{tt} -\Delta w_{tt} - \Delta w + V(x) w + w_{t}  = u_0 + u_1 - \Delta u_1
,\ \ \ (t,x)\in (0,\infty)\times {\bf R},\label{eq-w}\\
& w(0,x)=0, \quad  w_{t}(0,x)= u_{0}(x),\ \ \ x\in{\bf R}.\label{initial}\nonumber
\end{align}
Multiplying both sides of the equation in \eqref{eq-w} by $w_t$, and integrating over $[0,t]\times{\bf R}$ one has
\begin{align*}
&\dfrac{1}{2}\left(\|w_{t}(t,\cdot)\|^{2}+\|\nabla w_{t}(t,\cdot)\|^{2}+\|\nabla w(t,\cdot)\|^{2}+\|\sqrt{V(x)}w(t,\cdot)\|^{2}\right)
+ \int_{0}^{t}\|w_{s}(s,\cdot)\|^{2}\,ds\\
&= \dfrac{1}{2}\left(\|u_{0}\|^{2}+\|\nabla u_{0}\|^{2}\right) + \int_{0}^{t}\int_{\bf {R}}\left[u_{0}+u_{1}- \Delta u_{1}\right]w_{s}(s,x)\, dx\,ds\\
&= \dfrac{1}{2}\left(\|u_{0}\|^{2}+\|\nabla u_{0}\|^{2}\right)+\int_{0}^{t}\dfrac{d}{ds}\left(u_{0}+u_{1}- \Delta u_{1},\; w(s,\cdot) \right)_{L^{2}}ds\\
&= \dfrac{1}{2}\left(\|u_{0}\|^{2}+\|\nabla u_{0}\|^{2}\right)+\left(u_{0}+u_{1}-\Delta u_{1},\; w(t,\cdot) \right)_{L^{2}}\\
&\leq\dfrac{1}{2}\left(\|u_{0}\|^{2}+\|\nabla u_{0}\|^{2}\right)+\dfrac{1}{2}\left\|\dfrac{u_{0}+u_{1}-\Delta u_{1}}{\sqrt{V(x)}}\right\|^{2}+
\dfrac{1}{2}\|\sqrt{V(x)}w(t,\cdot)\|^{2}.
\end{align*}

The above identity implies the following inequality  
\[\dfrac{1}{2}\left(\|w_{t}(t,\cdot)\|^{2}+\|w_{xt}(t,\cdot)\|^{2}+\|w_{x}(t,\cdot)\|^{2}\right) + \int_{0}^{t}\|w_{s}(s,\cdot)\|^{2}\,ds\]
\[\leq \dfrac{1}{2}\left(\|u_{0}\|^{2}+\|\nabla u_{0}\|^{2} \right) +\dfrac{1}{2}\left\|\dfrac{u_{0}+u_{1}-\Delta u_{1}}{\sqrt{V(x)}}\right\|^{2}=J^{2}_{0},\]
which proves Lemma \ref{lem21} because of $w_{t}=u$.
\hfill \(\square\)
\vspace{0.1cm}\\
Next we prove the following intermediate result concerning the decay of the total energy. Conditions $\vert V'(x)\vert \leq \alpha V(x)$ and $\alpha^{2}\Vert V\Vert_{\infty} < 1$, as in Theorem \ref{theorem-1}, are unnecessary here just to obtain slower decay $E_{u}(t) = O(t^{-1})$ ($t \to \infty$).
\begin{pro}\label{pro21}\,Assume that $V \in {\rm BC}^{1}({\bf R})$ satisfies $V(x) > 0$ for $x \in {\bf R}$, and let $[u_{0},u_{1}] \in H^{2}({\bf R})\times H^{2}({\bf R})$ satisfying
$$\frac{u_{0} +u_{1}-\Delta u_1}{\sqrt{V(\cdot)}} \in L^2({\bf R}).$$
Then the strong solution $u \in C^{2}([0,\infty);H^{2}({\bf R}))$ to problem \eqref{eqn}-\eqref{initial}
satisfies the following asymptotic behavior 
\[E_{u}(t) \leq C(\Vert u_{0}\Vert_{H^{2}}^{2} + \Vert u_{1}\Vert_{H^{2}}^{2} + \left\|\dfrac{u_{0}+u_{1}-\Delta u_{1}}{\sqrt{V(\cdot)}}\right\|^{2})(1+t)^{-1},\quad (t \geq 0),\]
where $C = C(\Vert V\Vert_{\infty}, \Vert V'\Vert_{\infty}) > 0$ is a constant depending only on the quantities $\Vert V\Vert_{\infty}$ and $\Vert V'\Vert_{\infty}$.
\end{pro} 
\noindent
{\bf Proof:} To prove the proposition, we first note that it is standard to see the following basic estimate:
\[(1+t)E_{u}(t) \leq E_{u}(0) + \int_{0}^{t}\Vert u_{s}(s,\cdot)\Vert^{2}ds + \int_{0}^{t}\Vert \nabla u_{s}(s,\cdot)\Vert^{2}ds \]
\[+ \int_{0}^{t}\Vert\nabla u(s,\cdot)\Vert^{2}ds + \int_{0}^{t}\Vert \sqrt{V(\cdot)}u(s,\cdot)\Vert^{2}ds\]
\[\leq E_{u}(0) + \int_{0}^{t}\Vert u_{s}(s,\cdot)\Vert^{2}ds + \int_{0}^{t}\Vert \nabla u_{s}(s,\cdot)\Vert^{2}ds\]
\begin{equation}\label{1}
+ \int_{0}^{t}\Vert\nabla u(s,\cdot)\Vert^{2}ds + \Vert V\Vert_{\infty}\int_{0}^{t}\Vert u(s,\cdot)\Vert^{2}ds \quad (t \geq 0).
\end{equation}
Next, by multiplying both sides of \eqref{eqn} by $u(t,x)$ and integrating over ${\bf R}$ it follows that
\[\frac{d}{dt}(u_{t}(t,\cdot),u(t,\cdot)) + \Vert\nabla u(t,\cdot)\Vert^{2} + \Vert \sqrt{V(\cdot)}u(t,\cdot)\Vert^{2} + \frac{d}{dt}(\nabla u_{t}(t,\cdot),\nabla u(t,\cdot)) + \frac{1}{2}\frac{d}{dt}\Vert u(t,\cdot)\Vert^{2}\]
\[= \Vert u_{t}(t,\cdot)\Vert^{2} + \Vert \nabla u_{t}(t,\cdot)\Vert^{2}.\]
By integrating it on $[0,t]$, one has
\[(u_{t}(t,\cdot),u(t,\cdot)) + \int_{0}^{t}\Vert\nabla u(s,\cdot)\Vert^{2}ds + \int_{0}^{t}\Vert\sqrt{V(\cdot)}u(s,\cdot)\Vert^{2}ds + (\nabla u_{t}(t,\cdot),\nabla u(t,\cdot))+ \frac{1}{2}\Vert u(t,\cdot)\Vert^{2}\]
\begin{equation}\label{2}
= I_{0}^{2} + \int_{0}^{t}\Vert u_{s}(s,\cdot)\Vert^{2}ds + \int_{0}^{t}\Vert\nabla u_{s}(s,\cdot)\Vert^{2}ds,
\end{equation}
where
\[I_{0} := \left((u_{0},u_{1}) + (\nabla u_{0},\nabla u_{1}) + \frac{1}{2}\Vert u_{0}\Vert^{2}  \right)^{1/2}.\]
While, multiplying both sides of \eqref{eqn} by $-\Delta u_{t}$ and integrating over ${\bf R}\times [0,t]$ it follows that
\[\frac{1}{2}\Vert \nabla u_{t}(t,\cdot)\Vert^{2} + \frac{1}{2}\Vert\Delta u_{t}(t,\cdot)\Vert^{2}+ \frac{1}{2}\Vert\Delta u(t,\cdot)\Vert^{2}\]
\begin{equation}\label{3}
+\int_{0}^{t}\Vert\nabla u_{s}(s,\cdot)\Vert^{2}ds + \int_{0}^{t}\int_{{\bf R}}V(x)u(s,x)(-\Delta u_{s}(s,x))dxds = I_{1}^{2},
\end{equation}
where
\begin{equation}\label{I1}
I_{1} := \left(\frac{1}{2}\Vert\nabla u_{1}\Vert^{2} + \frac{1}{2}\Vert\Delta u_{1}\Vert^{2} + \frac{1}{2}\Vert\Delta u_{0}\Vert^{2}\right)^{1/2}.
\end{equation}
Note that
\[\int_{{\bf R}}V(x)u(s,x)(-\Delta u_{s}(s,x))dx = \int_{{\bf R}}u(s,x)\nabla u_{t}(s,x)\cdot\nabla V(x)dx + \frac{1}{2}\frac{d}{ds}\int_{{\bf R}}V(x)\vert\nabla u(s,x)\vert^{2}dx.\]
Then integration over $[0,\;t]$ implies that 
\[\int_{0}^{t}\int_{{\bf R}}V(x)u(s,x)(-\Delta u_{s}(s,x))dxds\]
\[= \int_{0}^{t}\int_{{\bf R}}u(s,x)\nabla u_{s}(s,x)\cdot\nabla V(x)dxds + \frac{1}{2}\int_{{\bf R}}V(x)\vert\nabla u(t,x)\vert^{2}dx - \frac{1}{2}\int_{{\bf R}}V(x)\vert\nabla u_{0}(x)\vert^{2}dx.\]
By substituting this equality into \eqref{3}, one can get the identity
\begin{align}\label{4}
&\frac{1}{2}\Vert \nabla u_{t}(t,\cdot)\Vert^{2} + \frac{1}{2}\Vert\Delta u_{t}(t,\cdot)\Vert^{2}+ \frac{1}{2}\Vert\Delta u(t,\cdot)\Vert^{2}+\int_{0}^{t}\Vert\nabla u_{s}(s,\cdot)\Vert^{2}ds\\
&+\int_{0}^{t}\int_{{\bf R}}u(s,x)\nabla u_{s}(s,x)\cdot\nabla V(x)dxds + \frac{1}{2}\int_{{\bf R}}V(x)\vert\nabla u(t,x)\vert^{2}dx = I_{1}^{2} + \frac{1}{2}\int_{{\bf R}}V(x)\vert\nabla u_{0}(x)\vert^{2}dx.\nonumber
\end{align}
Here, we need to use the estimate
\[\left\vert \int_{0}^{t}\int_{{\bf R}}u(s,x)\nabla u_{s}(s,x)\cdot\nabla V(x)dxds\right\vert \leq \int_{0}^{t}\int_{{\bf R}}\vert u(s,x)\vert\vert\nabla u_{s}(s,x)\vert\vert V'(x)\vert dxds\]
\[\leq \Vert V'\Vert_{\infty}\int_{0}^{t}\int_{{\bf R}}\vert u(s,x)\vert\vert\nabla u_{s}(s,x)\vert dxds\]
\begin{equation}\label{5}
\leq \varepsilon\Vert V'\Vert_{\infty}\int_{0}^{t}\Vert\nabla u_{s}(s,\cdot)\Vert^{2}ds + C_{\varepsilon}\Vert V'\Vert_{\infty}\int_{0}^{t}\Vert u(s,\cdot)\Vert^{2}ds
\end{equation}
with some constant $C_{\varepsilon} > 0$ which depends only on $\varepsilon > 0$. Thus, from \eqref{4} and \eqref{5} it follows that

\[\frac{1}{2}\Vert \nabla u_{t}(t,\cdot)\Vert^{2} + \frac{1}{2}\Vert\Delta u_{t}(t,\cdot)\Vert^{2}+ \frac{1}{2}\Vert\Delta u(t,\cdot)\Vert^{2}+(1-\varepsilon\Vert V'\Vert_{\infty})\int_{0}^{t}\Vert\nabla u_{s}(s,\cdot)\Vert^{2}ds\]
\begin{equation}\label{6}
+ \frac{1}{2}\int_{{\bf R}}V(x)\vert\nabla u(t,x)\vert^{2}dx \leq I_{1}^{2} + \frac{1}{2}\int_{{\bf R}}V(x)\vert\nabla u_{0}(x)\vert^{2}dx + C_{\varepsilon}\Vert V'\Vert_{\infty}\int_{0}^{t}\Vert u(s,\cdot)\Vert^{2}ds.
\end{equation}
On the other hand, it follows from \eqref{2} that
\[\int_{0}^{t}\Vert\nabla u(s,\cdot)\Vert^{2}ds + \int_{0}^{t}\Vert\sqrt{V(\cdot)}u(s,\cdot)\Vert^{2}ds + \frac{1}{2}\Vert u(t,\cdot)\Vert^{2}\]
\[
\leq I_{0}^{2} + \int_{0}^{t}\Vert u_{s}(s,\cdot)\Vert^{2}ds + \int_{0}^{t}\Vert\nabla u_{s}(s,\cdot)\Vert^{2}ds +\frac{1}{2}\Vert u_{t}(t,\cdot)\Vert^{2} + \frac{1}{2}\Vert u(t,\cdot)\Vert^{2} + \frac{1}{2}\Vert \nabla u_{t}(t,\cdot)\Vert^{2} + \frac{1}{2}\Vert \nabla u(t,\cdot)\Vert^{2},
\]
which implies
\begin{eqnarray}\label{7}
&&\int_{0}^{t}\Vert\nabla u(s,\cdot)\Vert^{2}ds + \int_{0}^{t}\Vert\sqrt{V(\cdot)}u(s,\cdot)\Vert^{2}ds\\
&\leq & I_{0}^{2} + \int_{0}^{t}\Vert u_{s}(s,\cdot)\Vert^{2}ds + \int_{0}^{t}\Vert\nabla u_{s}(s,\cdot)\Vert^{2}ds +\frac{1}{2}\Vert u_{t}(t,\cdot)\Vert^{2} + \frac{1}{2}\Vert \nabla u_{t}(t,\cdot)\Vert^{2} + \frac{1}{2}\Vert \nabla u(t,\cdot)\Vert^{2}.\nonumber
\end{eqnarray}
\noindent
Thus, for small $\varepsilon > 0$ , it follows from Lemma \ref{lem21}, \eqref{6}, \eqref{7} and \eqref{1-1} that 
\[(1-\varepsilon\Vert V'\Vert_{\infty})\int_{0}^{t}\Vert\nabla u_{s}(s,\cdot)\Vert^{2}ds\]
\begin{equation}\label{8}
\leq I_{1}^{2} + \frac{1}{2}\int_{{\bf R}}V(x)\vert\nabla u_{0}(x)\vert^{2}dx + C_{\varepsilon}\Vert V'\Vert_{\infty}J_{0}^{2} =: K_{0}^{2}
,
\end{equation}
and so
\begin{equation}\label{9}
\int_{0}^{t}\Vert\nabla u(s,\cdot)\Vert^{2}ds \leq I_{0}^{2} + E_{u}(0) + \frac{1}{1-\varepsilon\Vert V'\Vert_{\infty}}\left(I_{1}^{2} + \Vert\sqrt{V(\cdot)}\nabla u_{0}\Vert^{2}\right) + \frac{C_{\varepsilon}\Vert V'\Vert_{\infty}}{1-\varepsilon\Vert V'\Vert_{\infty}}J_{0}^{2} =: K_{1}^{2}.
\end{equation}

Finally, based on the estimates just obtained in Lemma \ref{lem21}, \eqref{1-1}, \eqref{1}, \eqref{8} and \eqref{9} the statement of Proposition \ref{pro21} can be easily obtained. Indeed, one can get
\begin{equation}\label{r-1}
(1+t)E_{u}(t) \leq 2E_{u}(0) + \Vert V\Vert_{\infty}J_{0}^{2} + \frac{K_{0}^{2}}{1-\varepsilon\Vert V'\Vert_{\infty}}K_{1}^{2} =: C_{1}^{2}.
\end{equation}

\begin{rem}{\rm The case of $\Vert V'\Vert_{\infty} = 0$ (that is $V(x) = V_{0} > 0$ for all $x \in {\bf R}$) is more easy to be treated, so it suffices to assume $\Vert V'\Vert_{\infty} \ne 0$ when checking the proof of Proposition \ref{pro21}.}
\end{rem}

\begin{rem}\label{rem22}
{\rm From \eqref{8}  in the proof of Proposition \ref{pro21} the following estimate holds,}
\begin{equation}\label{11}
\int_{0}^{t}\Vert\nabla u_{s}(s,\cdot)\Vert^{2}ds \leq \frac{K_{0}^{2}}{(1-\varepsilon\Vert V'\Vert_{\infty})}=:K_2^2.
\end{equation}
    \end{rem}

\section{ Fast Decay Estimates }
The idea to prove Theorem \ref{theorem-1} comes from the identity
$$\frac{d}{dt}\big[(1+t)^2E(t) \big]= 2(1+t)E(t) + (1+t)^2\frac{d}{dt}E(t), \;\; t>0,$$
which implies that 
\begin{align}\label{fast-ineq}
(1+t)^2E(t) \leq  E(0) + 2\int_0^t(1+s)E(s)ds, \;\; t>0,
\end{align}
because $E(t)$ is a non-increasing function of $t$.
\vspace{0.2cm}\\
Thus, to get the proof of faster decay of the total energy in Theorem \ref{theorem-1} it is necessary to estimate the integral in the right hand side of \eqref{fast-ineq}. That is, 
 from the definition of $E(t)$ one has to deal with the following inequality,
\[(1+t)^{2}E_{u}(t) \leq E_{u}(0) + \int_{0}^{t}(1+s)\Vert u_{s}(s,\cdot)\Vert^{2}ds + \int_{0}^{t}(1+s)\Vert \nabla u_{s}(s,\cdot)\Vert^{2}ds\]
\begin{equation}\label{31}
+ \int_{0}^{t}(1+s)\Vert\nabla u(s,\cdot)\Vert^{2}ds + \int_{0}^{t}(1+s)\Vert\sqrt{V(\cdot)}u(s,\cdot)\Vert^{2}ds.
\end{equation}
The first estimate one has to prove is the following.
\begin{lem}\label{lem31}
The strong solution satisfies the inequality
\begin{equation}\label{k3}
   \int_0^t (1+s)\Vert u_s(s,\cdot)\Vert^2ds \leq\; K_3^2
\end{equation}
where for some $0<\epsilon < 1$ determined in \eqref{r-1}:

$$K_3=\Big(2E(0) +  \Vert V\Vert_{\infty}J_{0}^{2} + \frac{K_{0}^{2}}{1-\varepsilon\Vert V'\Vert_{\infty}}K_{1}^{2} \Big)^{1/2}.$$

\end{lem}
{\bf Proof:}
By multiplying both sides of the equation \eqref{eqn} by $u_t$ it results that
$$\frac{1}{2}\frac{d}{dt}\Big( ||u_t(t,\cdot)||^2 + ||\nabla u_{t}(t,\cdot)||^2 + ||\nabla u(t,\cdot)||^2 + ||\sqrt{V(x)}u(t,\cdot)||^2\Big) + ||u_t(t,\cdot)||^2=  0 , \;\; t>0.$$
Now, multiplication this identity by $(1+t)$  implies that 
\[\frac{1}{2}\frac{d}{dt}\Big[(1+t)\big( 
||u_t(t,\cdot)||^2 + ||\nabla u_{t}(t,\cdot)||^2 + ||u_x(t,\cdot)||^2 + ||\sqrt{V(x)}u(t,\cdot)||^2\big)\Big] + (1+t)||u_t(t,\cdot)||^2 \]
\[ {\bf = }\frac{1}{2}\Big(||u_t(t,\cdot)||^{2} +||\nabla u_{t}(t,\cdot)||^2 + ||u_x(t,\cdot)||^2 + ||\sqrt{V(x)}u(t,\cdot)||^2 \Big), \;\; t>0.
\]
Integration the above identity over $[0,\;t]$ one gets for $t>0$
$$\frac{1}{2}\Big[(1+t)\big( 
||u_t(t,\cdot)||^2 + ||\nabla u_{t}(t,\cdot)||^2 + ||\nabla u(t,\cdot)||^2 + ||\sqrt{V(x)}u(t,\cdot)||^2\big)\Big] + \int_0^t(1+s)||u_s(s,\cdot)||^2ds$$ 
$$= E(0) + \int_0^t E(s)ds.$$
At this point, by re-examining the proof of Proposition \ref{pro21} we can see that it was already proven that
$$ \int_0^t E(s)ds \leq E(0) +  \Vert V\Vert_{\infty}J_{0}^{2} + \frac{K_{0}^{2}}{1-\varepsilon\Vert V'\Vert_{\infty}}K_{1}^{2}, \;\; t>0.$$
This estimate implies the proof of Lemma \ref{lem31}.
\hfill \(\square\)
\hspace{0.3cm}\\
To obtain the next estimate, the following identity can be derived by multiplying both sides of the equation \eqref{eqn} by $ (1+t)u$, integrating over $[0,t]\times {\bf R}$, and by the integration by parts, similarly to the proof of Proposition \ref{pro21}: 
\begin{align}\label{32}
&(1+t)(u_{t}(t,\cdot),u(t,\cdot)) - (u_{1}, u_{0}) = \frac{1}{2}\Vert u(t,\cdot)\Vert^{2} - \frac{1}{2}\Vert u_{0}\Vert^{2} + \int_{0}^{t}(1+s)\Vert u_{s}(s,\cdot)\Vert^{2}ds\nonumber\\
&-\int_{0}^{t}(1+s)\Vert\nabla u(s,\cdot)\Vert^{2}ds - \int_{0}^{t}(1+s)\Vert\sqrt{V(\cdot)}u(s,\cdot)\Vert^{2}ds-(1+t)(\nabla u_{t}(t,\cdot),\nabla u(t,\cdot)) \nonumber \\
&+ (\nabla u_{1},\nabla u_{0})
+ \frac{1}{2}\Vert \nabla u(t,\cdot)\Vert^{2} - \frac{1}{2}\Vert\nabla u_{0}\Vert^{2} +\int_{0}^{t}(1+s)\Vert \nabla u_{s}(s,\cdot)\Vert^{2}ds \nonumber \\
&-\frac{1}{2}(1+t)\Vert u(t,\cdot)\Vert^{2} + \frac{1}{2}\Vert u_{0}\Vert^{2} + \frac{1}{2}\int_{0}^{t}\Vert u(s,\cdot)\Vert^{2}ds.
\end{align}
By clean up this identity one has
\[\frac{1}{2}(1+t)\Vert u(t,\cdot)\Vert^{2} + \int_{0}^{t}(1+s)\Vert\nabla u(s,\cdot)\Vert^{2}ds + \int_{0}^{t}(1+s)\Vert\sqrt{V(\cdot)}u(s,\cdot)\Vert^{2}ds\]
\[= K_{0} - (1+t)(u_{t}(t,\cdot),u(t,\cdot)) + \int_{0}^{t}(1+s)\Vert u_{s}(s,\cdot)\Vert^{2}ds - (1+t)(\nabla u_{t}(t,\cdot),\nabla u(t,\cdot))\]
\begin{equation}\label{33}
+ \frac{1}{2}\Vert \nabla u(t,\cdot)\Vert^{2} + \int_{0}^{t}(1+s)\Vert \nabla u_{s}(s,\cdot)\Vert^{2}ds + \frac{1}{2}\int_{0}^{t}\Vert u(s,\cdot)\Vert^{2}ds + \frac{1}{2}\Vert u(t,\cdot)\Vert^{2},
\end{equation}
where
\[K_{0} := (\nabla u_{1},\nabla u_{0}) -\frac{1}{2}\Vert\nabla u_{0}\Vert^{2} + (u_{1},u_{0}).\]
Since
\[-(1+t)(u_{t}(t,\cdot),u(t,\cdot)) \leq (1+t)\Vert u_{t}(t,\cdot)\Vert^{2} + \frac{1+t}{4}\Vert u(t,\cdot)\Vert^{2},\]
and
\[- (1+t)(\nabla u_{t}(t,\cdot),\nabla u(t,\cdot)) \leq \frac{1+t}{2}\Vert\nabla u_{t}(t,\cdot)\Vert^{2} + \frac{1+t}{2}\Vert\nabla u(t,\cdot)\Vert^{2}\]
\[\leq (1+t)E_{u}(t) \leq C_{1}^{2}\quad(t \geq 0)\]
with some constant $C_{1} > 0$ (see \eqref{r-1}), it follows from \eqref{33} that
\[\frac{1}{4}(1+t)\Vert u(t,\cdot)\Vert^{2} + \int_{0}^{t}(1+s)\Vert\nabla u(s,\cdot)\Vert^{2}ds + \int_{0}^{t}(1+s)\Vert\sqrt{V(\cdot)}u(s,\cdot)\Vert^{2}ds\]
\begin{equation}\label{34}
\leq K_{0} + 3C_{1}^{2} + K(t) + \int_{0}^{t}(1+s)\Vert \nabla u_{s}(s,\cdot)\Vert^{2}ds + \int_{0}^{t}(1+s)\Vert u_{s}(s,\cdot)\Vert^{2}ds,
\end{equation}
where
\begin{equation}\label{Kt}
K(t) := \frac{1}{2}\Vert \nabla u(t,\cdot)\Vert^{2} + \frac{1}{2}\Vert u(t,\cdot)\Vert^{2} + \frac{1}{2}\int_{0}^{t}\Vert u(s,\cdot)\Vert^{2}ds.
\end{equation}
Here we denote the second energy by
\[2E^{*}(t) := \Vert \nabla u_{t}(t,\cdot)\Vert^{2} + \Vert \Delta u_{t}(t,\cdot)\Vert^{2} + \Vert \Delta u(t,\cdot)\Vert^{2}, \;\; t \geq 0.\]
Then, we have
\[\frac{d}{dt}E^{*}(t) + \big(V(\cdot)u(t,\cdot), -\Delta u_{t}(t,\cdot)\big) + \Vert \nabla u_{t}(t,\cdot)\Vert^{2} = 0.\]
Here, by integration by parts one has
\[(V(\cdot)u(t,\cdot), -\Delta u_{t}(t,\cdot)) = \int_{{\bf R}}u(t,x)\nabla V(x)\cdot\nabla u_{t}(t,x)dx + \frac{1}{2}\frac{d}{dt}\int_{{\bf R}}V(x)\vert\nabla u(t,x)\vert^{2}dx,\]
so that we get
\[\frac{d}{dt}E^{*}(t) + \int_{{\bf R}}u(t,x)\nabla V(x)\cdot\nabla u_{t}(t,x)dx + \frac{1}{2}\frac{d}{dt}\int_{{\bf R}}V(x)\vert\nabla u(t,x)\vert^{2}dx + \Vert \nabla u_{t}(t,\cdot)\Vert^{2} = 0.\]
Thus, it follows that
\[(1+t)E^{*}(t) + \frac{1+t}{2}\int_{{\bf R}}V(x)\vert\nabla u(t,x)\vert^{2}dx + \int_{0}^{t}(1+s)\Vert\nabla u_{s}(s,\cdot)\Vert^{2}ds\]
\begin{equation}\label{35}
= E^{*}(0) + \frac{1}{2}\int_{{\bf R}}V(x)\vert\nabla u_{0}(x)\vert^{2}dx + \int_{0}^{t}E^{*}(s)ds + \frac{1}{2}\int_{0}^{t}\int_{{\bf R}}V(x)\vert\nabla u(s,x)\vert^{2}dxds-F(t),
\end{equation}
where
\[F(t) := \int_{0}^{t}(1+s)\int_{{\bf R}}u(s,x)(\nabla u_{s}(s,x)\cdot\nabla V(x))dxds.\]

Let us estimate the function $F(t)$ by using the assumption on $V(x)$ such that
\[\vert V'(x)\vert \leq \alpha V(x)\]
with some $\alpha > 0$. Indeed,
\[\vert F(t)\vert \leq \int_{0}^{t}(1+s)\int_{{\bf R}}\vert u(s,x)\vert \;\vert\nabla u_{s}(s,x)\vert \;\vert V'(x)\vert dxds\]
\[\leq \int_{0}^{t}(1+s)\int_{{\bf R}}\vert u(s,x)\vert\vert\nabla u_{s}(s,x)\vert\alpha  V(x) dxds\]
\[\leq \frac{\alpha^{2}}{2}\int_{0}^{t}(1+s)\int_{{\bf R}}V(x)^{2}\vert u(s,x)\vert^{2} dxds + \frac{1}{2}\int_{0}^{t}(1+s)\int_{{\bf R}}\vert\nabla u_{s}(s,x)\vert^{2} dxds\]
\[\leq \frac{\alpha^{2}}{2}\Vert V\Vert_{\infty}\int_{0}^{t}(1+s)\int_{{\bf R}}V(x)\vert u(s,x)\vert^{2} dxds + \frac{1}{2}\int_{0}^{t}(1+s)\int_{{\bf R}}\vert\nabla u_{s}(s,x)\vert^{2} dxds,\]
which implies from \eqref{35} that
\[(1+t)E^{*}(t) + \frac{1+t}{2}\int_{{\bf R}}V(x)\vert\nabla u(t,x)\vert^{2}dx + \frac{1}{2}\int_{0}^{t}(1+s)\Vert\nabla u_{s}(s,\cdot)\Vert^{2}ds\]
\[\leq E^{*}(0) + \frac{1}{2}\int_{\bf R}V(x)\vert\nabla u_{0}(x)\vert^{2}dx + \int_{0}^{t}E^{*}(s)ds + \frac{1}{2}\int_{0}^{t}\int_{{\bf R}}V(x)\vert\nabla u(s,x)\vert^{2}dxds\]
\begin{equation}\label{36}
 + \frac{\alpha^{2}}{2}\Vert V\Vert_{\infty}\int_{0}^{t}(1+s)\int_{{\bf R}}V(x)\vert u(s,x)\vert^{2} dxds.
\end{equation}
From \eqref{36} one has
\[\int_{0}^{t}(1+s)\Vert\nabla u_{s}(s,\cdot)\Vert^{2}ds \leq 2E^{*}(0) + \int_{\bf R}V(x)\vert\nabla u_{0}(x)\vert^{2}dx + 2\int_{0}^{t}E^{*}(s)ds \]
\begin{equation}\label{37}
+ \int_{0}^{t}\int_{{\bf R}}V(x)\vert\nabla u(s,x)\vert^{2}dxds
 + \alpha^{2}\Vert V\Vert_{\infty}\int_{0}^{t}(1+s)\int_{{\bf R}}V(x)\vert u(s,x)\vert^{2} dxds.
\end{equation}
By combining \eqref{34} and \eqref{37} one can arrive at the meaningful bound:
\[\frac{1}{4}(1+t)\Vert u(t,\cdot)\Vert^{2} + \int_{0}^{t}(1+s)\Vert\nabla u(s,\cdot)\Vert^{2}ds + (1-\alpha^{2}\Vert V\Vert_{\infty})\int_{0}^{t}(1+s)\Vert\sqrt{V(\cdot)}u(s,\cdot)\Vert^{2}ds\]
\[\leq K_{0} + 3C_{1}^{2} + K(t) + 2E^{*}(0) + \int_{{\bf R}}V(x)\vert\nabla u_{0}(x)\vert^{2}dx + 2\int_{0}^{t}E^{*}(s)ds\]
\begin{equation}\label{38}
 + \int_{0}^{t}\int_{{\bf R}}V(x)\vert\nabla u(s,x)\vert^{2}dxds + \int_{0}^{t}(1+s)\Vert u_{s}(s,\cdot)\Vert^{2}ds.
\end{equation}
Thus, applying the estimates \eqref{9} and \eqref{k3} one has

\[\int_{0}^{t}(1+s)\Vert\nabla u(s,\cdot)\Vert^{2}ds + (1-\alpha^{2}\Vert V\Vert_{\infty})\int_{0}^{t}(1+s)\Vert\sqrt{V(\cdot)}u(s,\cdot)\Vert^{2}ds\]
\[\leq K_{0} + 3C_{1}^{2} + K(t) + 2E^{*}(0) + \int_{{\bf R}}V(x)\vert\nabla u_{0}(x)\vert^{2}dx + 2\int_{0}^{t}E^{*}(s)ds\]
\[+\,\,||V||_{\infty}\int_{0}^{t}\int_{{\bf R}}\vert\nabla u(s,x)\vert^{2}dxds + \int_{0}^{t}(1+s)\Vert u_{s}(s,\cdot)\Vert^{2}ds \]
\begin{equation}\label{mainbound}
\leq K_{0} + 3C_{1}^{2} + K(t) + 2E^{*}(0) + \int_{{\bf R}}V(x)\vert\nabla u_{0}(x)\vert^{2}dx + 2\int_{0}^{t}E^{*}(s)ds + ||V||_{\infty}K_1^2 + K_3^2.
\end{equation}
Because of Lemma \ref{lem21}, the definition of $K(t)$ in \eqref{Kt} implies that
$$K(t) \leq E(0) + \frac{1}{2}J_0^2 =:J_1^2. $$
Moreover, it is standard that $E^{*}(0)$ is bounded in terms of the norms of $u_0$ and $u_1$ in $H^2(\R)$.\\
Considering these two above estimates in \eqref{mainbound} one has arrived at the following results, which will help one to estimate the integrals on the right hand side of \eqref{31}.

\begin{lem}\label{39}\,
Under the assumption $\alpha^{2}\Vert V\Vert_{\infty} < 1$, it holds that
\[\int_{0}^{t}(1+s)\Vert\sqrt{V(\cdot)}u(s,\cdot)\Vert^{2}ds 
\leq C^{*}\Big( E_0^2 + 2\int_{0}^{t}E^{*}(s)ds \Big),\]
where $C^{*} := (1-\alpha^{2}\Vert V\Vert_{\infty})^{-1} > 0$ and 
$$E_0^2=:  K_{0} + 3C_{1}^{2} + J_1^2 + 2E^{*}(0) + \int_{{\bf R}}V(x)\vert\nabla u_{0}(x)\vert^{2}dx + ||V||_{\infty}K_1^2 + K_3^2. $$
Moreover,
\begin{align*}\label{mainbound2}
\int_{0}^{t}(1+s)\Vert\nabla u(s,\cdot)\Vert^{2}ds \leq E_0^2  + 2\int_{0}^{t}E^{*}(s)ds.
\end{align*}
\end{lem}
\vspace{0.2cm}
The next steps to complete the estimates in the above two lemmas is to obtain an upper bound to the integral 
\begin{equation}\label{final*-1}
2\int_0^tE^{*}(s)ds = \int_0^t ||\nabla u_s(s,\cdot)||^2ds + \int_0^t||\Delta u_s(s,\cdot)||^2ds + \int_0^t||\Delta u(s,\cdot)||^2ds.
\end{equation}\\
Note that the estimate to the first integral on the hight hand side of \eqref{final*-1} is already obtained by the estimate in \eqref{11}.\\
The estimates to the other two integrals will be done in the next lemmas but 
first we need to prove the following result for preparation.
\begin{lem}\label{est-Delta}
The $L^2$-norm of the  Laplacian of  $u(t,x)$ is bounded, that is
\begin{align}\label{lapla-est.2}
&\frac{1}{2}\Vert\Delta u(t,\cdot)\Vert^{2} 
\leq  I_{1}^{2} +  \frac{1}{2}\int_{{\bf R}}V(x)\vert\nabla u_{0}(x)\vert^{2}dx + C_{\varepsilon}\Vert V'\Vert_{\infty}J_0^2 =:J_2^2,
\end{align}
\end{lem}
where 
\[I_{1} := \left(\frac{1}{2}\Vert\nabla u_{1}\Vert^{2} + \frac{1}{2}\Vert\Delta u_{1}\Vert^{2} + \frac{1}{2}\Vert\Delta u_{0}\Vert^{2}\right)^{1/2}\]
is defined in \eqref{I1}.
\\
{\bf Proof:}
To prove the estimate to $\| \Delta u(t)\|^{2}$ we use the inequality in \eqref{6}. In fact, from that inequality one has  for some constant $C_{\varepsilon} > 0$ which depends only on $\varepsilon > 0$ such that 
$1-\varepsilon||V'||_{\infty}>0$, the  following estimate 
\begin{align}\label{66}
 \frac{1}{2}\Vert\Delta u(t,\cdot)\Vert^{2}
 \leq & I_{1}^{2} + \frac{1}{2}\int_{{\bf R}}V(x)\vert\nabla u_{0}(x)\vert^{2}dx + C_{\varepsilon}\Vert V'\Vert_{\infty}\int_{0}^{t}\Vert u(s,\cdot)\Vert^{2}ds\\
 \leq & I_{1}^{2} + \frac{1}{2}\int_{{\bf R}}V(x)\vert\nabla u_{0}(x)\vert^{2}dx + C_{\varepsilon}\Vert V'\Vert_{\infty}J_0^2,
\end{align}
by using  the estimate in Lemma \ref{lem21}.
\hfill \(\square\)\\
\vspace{0.2cm}
In the next lemma one needs the assumption on the $V'' \in {\rm BC}({\bf R})$ and must imposes further regularity on the initial data $(u_{0},u_{1}) \in H^{3}({\bf R})\times H^{3}({\bf R})$.
\begin{lem}\label{lem34} Let $u(t,x)$ be the regular solution of the Cauchy problem \eqref{eqn}--\eqref{initial}. Then,
\begin{equation*}
\dfrac{1}{2}\Vert \Delta u_{t}(t,\cdot)\Vert^{2} +\dfrac{1}{2}\int_{0}^{t}\|\Delta u_{s}(s,\cdot)\|^{2}\,ds \leq L^{2}_{0}, \quad t>0,
\end{equation*}
where\\
\vspace{-0.6cm}
$$L_{0}^{2} =: E_{2}(0) +C_{\delta}\|V''\|_{\infty}\, J^{2}_{0}+ 2C_{\delta}\|V'\|_{\infty}\,E(0)+
\|V\|_{\infty}\,J_{2}^{2} = L_{0}^{2},$$
and
\[2E_{2}(0) := \Big[ \Vert \Delta u_{1}\Vert^{2}+\Vert \Delta (u_{1})_{x}\Vert^{2}+\Vert \Delta (u_{0})_{x}\Vert^{2}\Big].\]

\end{lem}
{\bf Proof:}
Using the non-standard multiplier $\Delta^2 u_t$ for the wave equation, one has
\begin{align*}
\dfrac{1}{2}\left[\dfrac{d}{dt}\Vert \Delta u_{t}(t,\cdot)\Vert^{2}+\dfrac{d}{dt}\Vert \Delta u_{tx}(t,\cdot)\Vert^{2}
+\dfrac{d}{dt}\Vert \Delta u_{x}(t,\cdot)\Vert^{2} \right]+\|\Delta u_{t}(t,\cdot)\|^{2}\\
= -\Bigl(\Delta \bigr(V(x)u(t,\cdot)\bigr),\; \Delta u_{t}(t,\cdot)\Bigr)_{L^{2}}, \;\;t>0.
\end{align*}
Integration over $[0,\; t]$, $t>0$,  implies that for some $\delta > 0$
\begin{eqnarray*}
&&\dfrac{1}{2}\Big[\Vert\Delta u_{t}(t,\cdot)\Vert^{2} + \Vert\Delta u_{tx}(t,\cdot)\Vert^{2}+ \Vert \Delta u_{x}(t,\cdot)\Vert^{2}\Big]+
\int_{0}^{t}\|\Delta u_{s}(s,\cdot)\|^{2}\,ds\\
& \quad =&\dfrac{1}{2}\Big[ \Vert \Delta u_{1}\Vert^{2}+\Vert \Delta (u_{1})_{x}\Vert^{2}+\Vert \Delta (u_{0})_{x}\Vert^{2}\Big]
-\int_{0}^{t}\Bigl( \Delta \big(V(x)u(s,\cdot) \big) ,\; \Delta u_{s}(s,\cdot)\Bigl)_{L^{2}}\, ds\\
&\leq& E_{2}(0) -\int_{0}^{t}\Bigl(V''(x)u(s,\cdot)+2V'(x) \nabla u(s,\cdot)+V(x) \Delta u(s,\cdot),\; \Delta u_{s}(s,\cdot)\Bigr)_{L^{2}}\, ds\\
& \quad = &E_{2}(0) -\int_{0}^{t}\Bigl(V''(x)u(s,\cdot),\;\Delta u_{s}(s,\cdot)\Bigr)_{L^{2}}\, ds
 -\int_{0}^{t}\Bigl(2V'(x) \nabla u(s,\cdot),\; \Delta u_{s}(s,\cdot)\Bigr)_{L^{2}}\, ds\\
& - & \int_{0}^{t}\Bigl(V(x) \Delta u(s,\cdot),\; \Delta u_{s}(s,\cdot)\Bigr)_{L^{2}}\, ds\\
&\leq& E_{2}(0) +C_{\delta}\|V''\|^{2}\int_{0}^{t}\|u(s,\cdot)\|^{2}\,ds +\delta\int_{0}^{t}\|\Delta u_{s}(s,\cdot)\|^{2}\,ds
 + 2C_{\delta}\|V'\|_{\infty}\int_{0}^{t}\|\nabla u(s,\cdot)\|^{2}\,ds\\
& \quad+&\delta\int_{0}^{t}\|\Delta u_{s}(s,\cdot)\|^{2}\,ds - \int_{0}^{t}\Bigl(\sqrt{V(x)} \Delta u_{}(s,\cdot),\; \sqrt{V(x)} \Delta u_{s}(s,\cdot) \Bigr)\,ds,
\end{eqnarray*}
where $C_{\delta} > 0$ is a constant depending only on $\delta > 0$.\\
For simplification of symbols, $\Vert u(t)\Vert$ and $\Vert u\Vert$, etc., will be used in place of $\Vert u(t,\cdot)\Vert$ in the following. Then, taking in the last estimate $\delta >0 $ such that $1-2\delta>0$ one has for $t>0$,\\
\begin{align*}
&\dfrac{1}{2}\Big[\Vert \Delta u_{t}\Vert^{2} + \Vert\Delta u_{tx}\Vert^{2}+ \Vert \Delta u_{x}\Vert^{2}\Big]+
  (1-2\delta)\int_{0}^{t}\|\Delta u_{s}\|^{2}\,ds\\
&\leq E_{2}(0) +C_{\delta}\|V''\|_{\infty}\int_{0}^{t}\|u\|^{2}\,ds+
2C_{\delta}\|V'\|_{\infty}\int_{0}^{t}\|\nabla u\|^{2}\,ds -
\dfrac{1}{2}\int_{0}^{t}\dfrac{d}{ds}\Big(\left\|\sqrt{V(\cdot)} \Delta u(s)\right\| \Big)\,ds\\
&  =E_{2}(0) +C_{\delta}\|V''\|_{\infty}\int_{0}^{t}\|u\|^{2}\,ds+
2C_{\delta}\|V'\|_{\infty}\int_{0}^{t}\|\nabla u\|^{2}\,ds -
\dfrac{1}{2}\left\|\sqrt{V(\cdot)} \Delta u\right\|^{2}+ \dfrac{1}{2}\left\|\sqrt{V(\cdot)}\Delta u_{0}\right\| \\
&\leq E_{2}(0) +C_{\delta}\|V''\|_{\infty}\, J^{2}_{0}+ 2C_{\delta}\|V'\|_{\infty}\,K_{1}^{2}
-\dfrac{1}{2}\left\|\sqrt{V(\cdot)} \Delta u\right\|^{2} +\dfrac{1}{2}\left\|\sqrt{V(\cdot)} \Delta u_0\right\|^{2} \\
&\leq E_{2}(0) +C_{\delta}\|V''\|_{\infty}\, J^{2}_{0}+ 2C_{\delta}\|V'\|_{\infty}\,K_{1}^{2}+\frac{1}{2}\Vert \sqrt{V(\cdot)}\Delta u_{0}\Vert^{2} =: L_{0}^{2},
\end{align*}
due to Lemma \ref{lem21} and estimate \eqref{9}.
\vspace{0.2cm}\\
At this point, taking $ \delta = \frac{1}{4}$ one obtain the proof of  above lemma. Note that with $L_0$ is  a constant independent of the  time $t$.
Incidentally, the calculations so far can actually be justified by the density arguments.
\hfill \(\square\)
\vspace{0.3cm}\\
Finally one has to prove the following lemma.
\begin{lem}\label{lem:2-4}
 The regular solution $u(t,x)$ of the problem \eqref{eqn}--\eqref{initial} satisfies
\begin{equation*}
\int_{0}^{t}\|\Delta u(s,\cdot)\|^{2}\,ds \leq L^{2}_{1}
\end{equation*}
with $L^{2}_{1}:= C_0+ E(0) + K_2^2+ I_1^2 +  \dfrac{1}{2}\|V||_{\infty}^2 \big( 2J_0^2 + 2L_0^2 \big) + L_0^2. $ The constants in the definition of $L_1$ are defined in this and  previous sections,
where the newly appeared constant $C_{0} > 0$ can be defined by
\[C_{0} := \big(\nabla u_{1},\; \nabla u_{0}\big)+\big(\Delta u_{1},\; \Delta u_{0}\big)
+ \dfrac{1}{2}\|\nabla u_{0}\|^{2}.\]
\end{lem}
{\bf Proof:}
Multiplying both sides of \eqref{eqn} by $-\Delta u$ and integrating over ${\bf R}$ one has
\begin{align*}
&\dfrac{d}{dt}\Big[\left( \nabla u_{t},\; \nabla u\right)+\left( \Delta u_{t},\; \Delta u\right) \Big]
-\|\nabla u_{t}\|^{2}-\| \Delta u_{t}\|^{2}+\|\Delta u\|^{2}-\left(V(x)u,\;\Delta u\right)\\
&+\dfrac{1}{2}\dfrac{d}{dt}\|\nabla u\Vert^{2}=0, \quad t>0.
\end{align*}
\vspace{0.2cm}
Then, by integrating on $[0,\;t]$ one has
\begin{align*}
&\left( \nabla u_{t},\; \nabla u\right)+\left(\Delta u_{t},\;\Delta u\right)-\int_{0}^{t}\Big[
\|\nabla u_{s}(s,\cdot)\|^{2}+\|\Delta u_{s}(s,\cdot)\|^{2} \Big]\, ds+\int_{0}^{t}\|\Delta u\|^{2}\,ds
+\dfrac{1}{2}\|\nabla u\|^{2}\\
&=\dfrac{1}{2}\|\nabla u_{0}\|^{2}+\left( \nabla u_{1},\; \nabla u_{0}\right)
+\left(\Delta u_{1},\; \Delta u_{0}\right)+\int_{0}^{t}\left(V(x)u,\; \Delta u\right)\,ds.
\end{align*}
Then, rearranging the terms of above identity and using \eqref{1-1}, one has
\begin{align*}
\dfrac{1}{2}\| \nabla u\|^{2}+\int_{0}^{t}\| \Delta u\|^{2}\,ds 
&= 
\big(\nabla u_{1},\; \nabla u_{0}\big)+\big(\Delta u_{1},\; \Delta u_{0}\big)
+ \dfrac{1}{2}\|\nabla u_{0}\|^{2}
-\left( \nabla u_{t},\; \nabla u\right)-\left(\Delta u_{t},\;\Delta u\right)\\
& \quad  + \int_{0}^{t}\Big[\|\nabla u_{s}(s,\cdot)\|^{2}+\|\Delta u_{s}(s,\cdot)\|^{2}\Big]\,ds+ \int_{0}^{t} \big( V(x)u, \Delta u\big)\,ds\\
&\leq
C_{0}+ \dfrac{1}{2}\|\nabla u_{t}\|^{2}+\dfrac{1}{2}\|\nabla u\|^{2}
+\dfrac{1}{2}\| \Delta u_{t}\|^{2}+\dfrac{1}{2}\|\Delta u\|^{2} +\int_{0}^{t}\| \nabla u_{s}(s,\cdot)\|^{2}\,ds\\
& \quad  +\int_{0}^{t}\|\Delta u_{s}(s,\cdot)\|^{2}\,ds +\dfrac{1}{2}\int_{0}^{t}\|V(x)u\|^{2}\,ds +\dfrac{1}{2}\int_{0}^{t}\|\Delta u\|^{2}\,ds\\
&\leq C_{0}+E(0)+\dfrac{1}{2}\|\Delta u_{t}\|^{2}+\dfrac{1}{2}\|\Delta u\|^{2} +\int_{0}^{t}\|\nabla u_{s}(s,\cdot)\|^{2}\,ds \\
& \quad +\int_{0}^{t}\|\Delta u_{s}(s,\cdot)\|^{2}\,ds +\dfrac{1}{2}\|V||_{\infty}^2\int_{0}^{t}\|u\|^{2}\,ds +\dfrac{1}{2}\int_{0}^{t}\|\Delta u\|^{2}\,ds.
\end{align*}
\noindent
This inequality implies the following estimate for $t>0$
\begin{align}\label{lapla-est}
&\dfrac{1}{2}\|\nabla u\|^{2}+ \dfrac{1}{2}\int_{0}^{t}\| \Delta u\|^{2}\,ds 
\leq C_{0}+E(0)+\dfrac{1}{2}\|\Delta u_{t}\|^{2}+\dfrac{1}{2}\|\Delta u\|^{2} \nonumber\\
& \quad +\int_{0}^{t}\|\nabla u_{s}(s,\cdot)\|^{2}\,ds+\int_{0}^{t}\|\Delta u_{s}(s,\cdot)\|^{2}\,ds +\dfrac{1}{2}\|V||_{\infty}^2\int_{0}^{t}\|u\|^{2}\,ds \nonumber  \\
&\leq C_0 + E(0) + L_0^2 + K_2^2 +  \dfrac{1}{2}\|V||_{\infty}^2J_0^2 + \dfrac{1}{2}\| \Delta u\|^{2}.
\end{align}
by Lemma~\ref{lem34}, estimate \eqref{11} and  Lemma~\ref{lem21}.\\
To get the upper bound estimate to $\| \Delta u(t)\|^{2}$ we use the inequality in \eqref{3}. In fact, from that inequality we have
\begin{align}\label{lapla-est.2}
&\frac{1}{2}\Vert\Delta u(t,\cdot)\Vert^{2} 
\leq I_{1}^{2} -\int_{0}^{t}\int_{{\bf R}}V(x)u(s,x)(-\Delta u_{s}(s,x))dxds\nonumber\\
& \leq  I_{1}^{2}  + \frac{1}{2}\|V||_{\infty} \int_{0}^{t}\|u(s)\|^2\;ds 
+ \frac{1}{2}\|V\|_{\infty} \int_{0}^{t}\;\Vert\Delta u_{s}(s,\cdot)\|^2ds  \nonumber \\
&\leq  I_{1}^{2} + \frac{1}{2}\|V\|_{\infty} \big( J_0^2 +  2L_0^2
\big),
\end{align}
by Lemmas \ref{lem21} and \ref{lem34}, where the constant $I_1$ depending only on the initial data is again given by
\[I_{1} := \left(\frac{1}{2}\Vert\nabla u_{1}\Vert^{2} + \frac{1}{2}\Vert\Delta u_{1}\Vert^{2} + \frac{1}{2}\Vert\Delta u_{0}\Vert^{2}\right)^{1/2}.\]
Now, combining estimates \eqref{lapla-est}   and \eqref{lapla-est.2} the proof of Lemma \ref{lem:2-4} follows
with
\begin{align*}
L_1^2 &=C_0 + E(0) + L_0^2 + K_2^2 +  \dfrac{1}{2}\|V||_{\infty}^2J_0^2 + I_{1}^{2} + \frac{1}{2}\|V\|_{\infty} \big( J_0^2 +  2L_0^2 \big)\\
&= C_0+ E(0) + K_2^2+ I_1^2 +  \dfrac{1}{2}\|V||_{\infty}^2 \big( 2J_0^2 + 2L_0^2 \big) + L_0^2.
\end{align*}. 
 \hfill \(\square\)

Let us finalize the proof of Theorem \ref{theorem-1}.\\
{\it Proof of Theorem \ref{theorem-1} completed.}\\
Now, using the estimates in Remark 2.2 and  Lemmas \ref{lem34} and \ref{lem:2-4} we get the important estimate 
\begin{equation}\label{final*}
\int_0^tE^{*}(s)ds \leq K_2^2 + L_0^2 + L_1^2 =:C_2^2,
\end{equation}
and from Lemma \ref{39} the following estimates
\begin{equation}\label{100}
\int_{0}^{t}(1+s)\Vert\sqrt{V(\cdot)}u(s,\cdot)\Vert^{2}ds 
\leq C^{*}\Big( E_0^2 + 2 C_2^2 \Big)
\end{equation}
and
\begin{align}\label{mainbound2}
\int_{0}^{t}(1+s)\Vert\nabla u(s,\cdot)\Vert^{2}ds \leq E_0^2  + 2C_2^2.
\end{align}
can be obtained under the assumption that $\alpha^{2}\Vert V\Vert_{\infty} < 1$.\\
\vspace{0.2cm}
To get the fast decay for the energy, using the estimate \eqref{k3}  in Lemma \ref{lem31} and estimates \eqref{100} and \eqref{mainbound2} we conclude from \eqref{31}

\[(1+t)^{2}E_{u}(t) \leq E_{u}(0) + K_3^2  + \int_{0}^{t}(1+s)\Vert \nabla u_{s}(s,\cdot)\Vert^{2}ds\]
\begin{equation}\label{101}
+  E_0^2  + 2C_2^2 +C^{*}\Big( E_0^2 + 2 C_2^2 \Big).
\end{equation}
Finally, to conclude the proof of the fast  decay for the total energy we get an upper bound to the integral on the right hand side of \eqref{101}.
In fact, from \eqref{37}, \eqref{final*}, \eqref{9} and \eqref{100}, one can get the desired estimate to the integral on the right hand side of \eqref{101}:
\begin{align*}
&\int_{0}^{t}(1+s)\Vert\nabla u_{s}(s,\cdot)\Vert^{2}ds\\
&\leq 2E^{*}(0) + \int_{\bf R}V(x)\vert\nabla u_{0}(x)\vert^{2}dx + 2C_2^2  
+ 2||V||_{\infty}\int_{0}^{t}\int_{{\bf R}}\vert\nabla u(s,x)\vert^{2}dxds\\
& + \alpha^{2}\Vert V\Vert_{\infty}\int_{0}^{t}(1+s)\int_{{\bf R}}V(x)\vert u(s,x)\vert^{2} dxds\\
&\leq 2E^{*}(0) + \int_{\bf R}V(x)\vert\nabla u_{0}(x)\vert^{2}dx + 2C_2^2  
+ 2||V||_{\infty}K_1^2 +  \alpha^{2}\Vert V\Vert_{\infty} C^{*}\Big( E_0^2 + 2 C_2^2 \Big).
\end{align*}
Incidentally, in order to get the $L^{2}$-decay estimate in Theorem \ref{theorem-1}, it suffices to combine the several preliminary estimates derived in \eqref{1-1}, \eqref{k3}, \eqref{34}, Lemma \ref{lem21} and the inequality just obtained above.
\hfill
$\Box$

\begin{rem}\,{\rm In order to get slower decay $E_{u}(t) = O(t^{-1})$ ($t \to \infty$), one has assumed the regularity $u_{j} \in H^{2}({\bf R})$ ($j = 0,1$), and to get faster decay $E_{u}(t) = O(t^{-2})$ ($t \to \infty$), one has just imposed the assumption $u_{j} \in H^{3}({\bf R})$ ($j = 0,1$). This may reflect a kind of regularity loss structure of the equation with the rotational inertia term $-\Delta u_{tt}$. }
\end{rem} 


\section{Appendix}

In this appendix we discuss the existence of regular solution and its regularity issue to problem \eqref{eqn}-\eqref{initial}. We apply the Lumer-Phillips theorem via semi-group theory. For generalization, we treat the problem \eqref{eqn}-\eqref{initial} in the $n$-dimensional space ${\bf R}^{n}$.\\

Define the operator ${\cal A}: D({\cal A}) = {\cal H} \to {\cal H} $ by
\begin{gather*}
	\mathcal{A} =	\begin{pmatrix}
		0& I \\ -J(I+A)& 0
	\end{pmatrix},~~~
\end{gather*}
where $A := -\Delta$ with $D(A) = H^{2}({\bf R}^{n})$ with the identity operator $I$ in $L^{2}({\bf R}^{n})$, and we set $J := (I+A)^{-1}$, ${\cal H} := H^{2}({\bf R}^{n})\times H^{2}({\bf R}^{n})$. Note that $J(I+A)u = u$ for $u \in D(A)$.

Setting 
\begin{gather*}
	U(t) =	\begin{pmatrix}
		u(t)\\v(t)
	\end{pmatrix},~~~
	U_0=	\begin{pmatrix}
		u_0\\u_1
	\end{pmatrix},~~~
\end{gather*}
the problem \eqref{eqn}-\eqref{initial} can be converted into the following Cauchy one in ${\cal H}$:
\begin{equation}\label{A1}
U_{t}(t) = {\cal A}U(t) + L_{V}U(t) + FU(t), t > 0, 
\end{equation}
and
\begin{equation}\label{A2}
U(0) = U_{0}.
\end{equation}
Here, we have just defined two operators $L_{V}:{\cal H} \to {\cal H}$ and $F:{\cal H} \to {\cal H}$ as follows:
\[L_{V}\begin{pmatrix}
		u\\v
	\end{pmatrix}: = \begin{pmatrix}
		0\\JV(\cdot)u
	\end{pmatrix},\quad F\begin{pmatrix}
		u\\v
	\end{pmatrix} := \begin{pmatrix}
		0\\-J(v-u)		
	\end{pmatrix}\] 
for each $\begin{pmatrix}
		u\\v
	\end{pmatrix} \in {\cal H}$. For this expression it should be noticed the following equivalent expression to the equation \eqref{eqn}:
\[u_{tt} + J(I+A)u + JV(x)u + J(u_{t}-u) = 0.\]	
Note that the mapping $L_{V}$ is well-defined because $V(\cdot)u \in L^{2}({\bf R}^{n})$ for $u \in H^{2}({\bf R}^{n})$ by the assumption $V \in L^{\infty}({\bf R}^{n})$.\\

Let us check the dissipativity of the operator ${\cal A}$. In fact, one can see that for each $\begin{pmatrix}u\\v\end{pmatrix} \in {\cal H}$, 
\[<{\cal A}U,U>_{{\cal H}} = <\begin{pmatrix}v\\-u\end{pmatrix}, \begin{pmatrix}u\\v\end{pmatrix}>_{{\cal H}} = (v,u)_{H^{2}}-(u,v)_{H^{2}} = 0,\]
where $<U,V>_{{\cal H}}$ and $(u,v)_{H^{2}}$ mean the standard inner product in ${\cal H}$ and $H^{2}$ spaces respectively.

Next, we can check the surjectivity of the mapping ${\cal I}-{\cal A}: D({\cal A}) \to {\cal H}$, where ${\cal I}$ is the identity on ${\cal H}$, that is, we should find a vector $U := \begin{pmatrix}u\\v\end{pmatrix} \in D({\cal A})$ for each $G:= \begin{pmatrix}f\\g\end{pmatrix} \in {\cal H}$ such that
\[U-{\cal A}U = G.\]
But, this is quite easy because of the property $J(I+A)u = u$ for $u \in D(A)$. We shall omit the details. Therefore, the operator ${\cal A}: D({\cal A}) = {\cal H} \to {\cal H}$ is m-dissipative. By the Lumer-Phillips Theorem it is a generator of a $C_{0}$-semigroup of contractions on ${\cal H}$ (see \cite[Section 1.4, Theorem 4.3]{pazy}). Furthermore, by a perturbation theory of the $C_{0}$-semigroup (see \cite[Section 3.1, Theorem 1.1]{pazy}) the operator ${\cal A} + L_{V} + F$ is a generator of a $C_{0}$-semigroup. Note that the operator $L_{V} + F$ is bounded in ${\cal H}$, and this fact can be easily checked as follows. Indeed, we first note that the $H^{2}$-norm is equivalent to the graph norm of the Laplacian $A$, that is, $\Vert u\Vert_{H^{2}} \sim \Vert u\Vert + \Vert Au\Vert$ for each $u \in H^{2}({\bf R}^{n})$. Furthermore, we use two facts about the Yosida approximation such that 
\[AJw = w-Jw, \quad \Vert Jw\Vert_{L^{2}} \leq \Vert w\Vert_{L^{2}}\]
for $w \in L^{2}({\bf R}^{n})$. Then, for $U := \begin{pmatrix}u\\v\end{pmatrix} \in {\cal H}$, we see that
\[\Vert L_{V}U\Vert_{{\cal H}} = \Vert JV(\cdot)u\Vert_{H^{2}} \sim \Vert JV(\cdot)u\Vert_{L^{2}} + \Vert AJV(\cdot)u\Vert_{L^{2}}\]
\[= \Vert JV(\cdot)u\Vert_{L^{2}} + \Vert V(\cdot)u-JV(\cdot)u\Vert_{L^{2}} \leq C\left(\Vert JV(\cdot)u\Vert_{L^{2}} + \Vert V(\cdot)u\Vert_{L^{2}}\right)\]
\[\leq C'\Vert V(\cdot)u\Vert_{L^{2}} \leq C'\Vert V\Vert_{L^{\infty}}\Vert u\Vert_{L^{2}} \leq C'\Vert V\Vert_{L^{\infty}}\Vert U\Vert_{{\cal H}},\]
which implies the boundedness of the operator $L_{V}$ in ${\cal H}$. Similarly, one can easily check the boundedness of the operator $F$ in ${\cal H}$.  Thus, for each initial data $U_{0} = \begin{pmatrix}u_{0}\\u_{1}\end{pmatrix} \in {\cal H} = H^{2}({\bf R}^{n})\times H^{2}({\bf 
R}^{n})$, 
there exists a unique solution $U(t) = \begin{pmatrix}u(t)\\v(t)\end{pmatrix} \in C([0,\infty);D({\cal A}))\cap C^{1}([0,\infty);{\cal H}) = C^{1}([0,\infty);{\cal H})$ to problem \eqref{A1}-\eqref{A2}. This implies $v=u_t$ and $u \in C^{2}([0,\infty);H^{2}({\bf R}^{n}))$ for the solution $u(t,x)$ of the problem \eqref{eqn}-\eqref{initial}.

\begin{rem}{\rm By a similar consideration to above discussion, if we rechoose ${\cal H} := H^{4}({\bf R}^{4})\times H^{4}({\bf R}^{4}) = D({\cal A})$, one can obtain $u \in C^{2}([0,\infty);H^{4}({\bf R}^{n}))$ for the solution $u(t,x)$ of the problem \eqref{eqn}-\eqref{initial}, which is sufficiently useful in our discussion.}
\end{rem} 

\noindent{\em Acknowledgement.}
\smallskip The work of the first author (R. C. Char\~ao) was partially supported by Print/Capes - Process 88881.310536/2018-00 and the work of the second author (R. Ikehata) was supported in part by Grant-in-Aid for Scientific Research (C)20K03682  of JSPS. 



\end{document}